\documentstyle[12pt,amssymb]{amsart}

\title[Finite Quantum Groups]
{\Large \bf Rieffel Type Discrete Deformation of Finite Quantum Groups}
\author[Shuzhou Wang]
{\bf Shuzhou Wang}
\address{Department of Mathematics, University of California,
Berkeley, CA 94720
\newline \indent
Fax: 510-642-8204
}
\email{szwang@@math.berkeley.edu}

\date{}

\newtheorem{DF}{Definition}[section]
\newtheorem{LM}[DF]{Lemma}
\newtheorem{PROP}[DF]{Proposition}
\newtheorem{TH}[DF]{Theorem}
\newtheorem{COR}[DF]{Corollary}
\newtheorem{RMK}[DF]{Remark}
\newtheorem{RMKS}[DF]{Remarks}
\newtheorem{PROB}[DF]{Problem}

\newcommand{\bgdf}{\begin{DF}}
\newcommand{\nddf}{\end{DF}}
\newcommand{\bglm}{\begin{LM}}
\newcommand{\ndlm}{\end{LM}}
\newcommand{\bgprop}{\begin{PROP}}
\newcommand{\ndprop}{\end{PROP}}
\newcommand{\bgth}{\begin{TH}}
\newcommand{\ndth}{\end{TH}}
\newcommand{\bgcor}{\begin{COR}}
\newcommand{\ndcor}{\end{COR}}
\newcommand{\bgrmk}{\begin{RMK}}
\newcommand{\ndrmk}{\end{RMK}}
\newcommand{\bgrmks}{\begin{RMKS}}
\newcommand{\ndrmks}{\end{RMKS}}
\newcommand{\bgprob}{\begin{PROB}}
\newcommand{\ndprob}{\end{PROB}}
\newcommand{\bgeq}{\begin{eqnarray}}
\newcommand{\ndeq}{\end{eqnarray}}
\newcommand{\bgeqq}{\begin{eqnarray*}}
\newcommand{\ndeqq}{\end{eqnarray*}}

\newcommand{\QED}{\hfill Q.E.D.} 
\newcommand{\vv}{\vspace{4mm}\\}

\def\RR{{\Bbb R}}

\def\ZZ{{\Bbb Z}}

\def\JJ{{\hbar J}}

\numberwithin{equation}{section}

\newcommand{\dfref}[1]{Definition~\ref{#1}}

\newcommand{\propref}[1]{Proposition~\ref{#1}}
\newcommand{\thmref}[1]{Theorem~\ref{#1}}
\newcommand{\thref}[1]{Theorem~\ref{#1}}

\newcommand{\probref}[1]{Problem~\ref{#1}}

\begin{document}

\begin{abstract}
We introduce a discrete deformation of Rieffel type for finite (quantum)
groups. Using this, we give an example of a finite quantum group
$A$ of order $18$ such that
neither $A$ nor its dual can be expressed as a crossed product of the form
$C(G_1) \rtimes_{\tau} G_2$ with $G_1$ and $G_2$ ordinary finite groups.
We also give a deformation of finite groups of Lie type
by using their maximal abelian subgroups.
\end{abstract}

\maketitle

\section{
Introduction}

Since the work of Woronowicz \cite{Wor4,Wor5}, the theory of compact quantum
groups, notably the deformation theory of compact Lie groups, has been
intensively studied and is now quite well understood
(see e.g. \cite{Wor6,VS,LS1,Rf6,W1,W5}). However, this is
not the case for finite quantum groups.
Both as objects of great mathematical interest, like finite
groups, and as objects with potential important applications
in theoretical physics \cite{Cn,Cn3}, the theory of finite quantum groups
calls for more efforts of study.
To start with, the theory needs an interesting supply of examples,
which are still lacking so far, though a few non-trivial
examples have been studied \cite{KP2,BS2,EnockVain1,Nik}.

In this paper, we construct a class of finite quantum groups
by introducing a discrete deformation of Rieffel type for finite (nonabelian)
groups. In fact, just as in our earlier paper \cite{W4}, this
deformation can be applied to finite quantum groups as well, not just finite
groups. This construction is motivated by Rieffel's deformation of
compact Lie groups \cite{Rf8}, which has its origins in the Weyl-von Neumann
quantization (also called Moyal product) (cf \cite{Rf6}).
As a matter of fact, our formula for the discretely deformed product
(see \dfref{deformedproduct}) is an exact analog of the
product formula of von Neumann and Rieffel \cite{Neumann1,Rf6}.
In \cite{Rf8} (resp. \cite{W4}), actions of finite dimensional vector
spaces are used to deform Lie groups (resp. compact quantum groups) into new
quantum groups. In this paper, we use actions of finite abelian groups
to deform finite groups (and finite quantum groups) into new finite quantum
groups. Because of the nature of the objects we deal with,
we are spared the analytical complications met in \cite{Rf6,Rf8,W4} for the
actions of continuous abelian groups (viz. vector spaces).
Hence the arguments in this paper are of a purely algebraic nature.
Though the constructions of this paper are direct analogs of
\cite{Rf6,Rf8,W4} adapted to the actions of finite abelian groups,
the proofs of the main results given there do not directly generalize to the
new situation. Thus we have to develop different proofs for our main results.
The main cause of this is that many facts on Euclidean geometry
of $\RR^d$ as used in \cite{Rf6} do not have generalizations to
finite abelian groups (e.g. orthogonal complements, polar decomposition of
operators, etc), though one can develop to
some extend the ``Euclidean geometry''
on a finite abelian group with ``inner product'' given by a pairing which
identifies itself with its Pontryagin dual.

The construction of deformation in this paper {\em in the dual form}
is an example of Drinfeld's twistings \cite{Dr6},
just as the constructions in \cite{Rf8,Rf9,W4} (cf also \cite{Lev1,LS1}).
This again shows the relationship between the Rieffel type
deformations (generalizations of the Weyl-von Neumann quantization) and
the Drinfeld's twistings. For Kac algebras, the most general form
of twistings in the sense of Drinfeld \cite{Dr6} is studied by
Enock and Vainerman \cite{EnockVain1,Vainerman2},
following the work of Landstad
and Raeburn \cite{LandstadRaeburn1,Landstad1} on deformations
of locally compact groups.
Hence in the dual picture our construction constitutes a {\em distinguished
class} of twistings of Kac algebras in the sense of Enock and Vainerman
\cite{EnockVain1,Vainerman2}.
Instead of imposing rather complicated cocycle conditions
in addition to the existence of an abelian subgroup, such as the approach
in \cite{Landstad1,EnockVain1,Vainerman2}, our construction of deformation
is {\em canonically} associated with the abelian subgroup,
and it is a natural generalization of the Weyl-von Neumann-Rieffel
deformation. As a matter of fact, the twist $F$
(see formula \eqref{pseudo-2-cocycle}) for the dual of our construction
does not satisfy the 2-cocycle condition on the Kac algebra, but the
pseudo-2-cocycle condition, which is equivalent to the condition
that the associated twisted coproduct is coassociative, a minimal
requirement. It is interesting to note that it is not clear how to
see that $F$ is a pseudo-2-cocycle directly in the dual picture!
Also, unlike \cite{BS2,EnockVain1}, our construction does
not give rise to the $8$ dimensional quantum group of
Kac-Palyutkin \cite{KP2}. Note also that in \cite{Nik},
more specific examples along the lines of
\cite{Landstad1,EnockVain1,Vainerman2} are
given; it is also shown there that the $K_0$ ring of Hopf algebra is invariant
under twists, which is obvious for our construction.

The plan of this paper is as follows. As preparation for Sect. 3,
we construct in Sect. 2 a deformed $C^*$-algebra
$A_J$ for every finite dimensional unital $C^*$-algebra
$A$ that is endowed with an action $\alpha$ of a finite abelian
group $H$, where $J$ is a skew-symmetric automorphism on $H$.
See \thref{theoremalgebra}.
The construction in this section parallels the one in \cite{Rf6}.
In Sect. 3, for every finite quantum group $A$ containing a
finite abelian subgroup $T$, we construct an action $\alpha$ of
$H$ on the $C^*$-algebra $A$ and show that the deformation $A_J$ is also
a finite quantum group containing $T$ as a subgroup,
where $H = T \oplus T$, $J = S \oplus (-S)$, and $S$ is a skew-symmetric
automorphism on $T$. See the main result \thref{theoremquantumgroup}.
This theorem parallels the main results in \cite{Rf8,W4},
and is announced in Section 2 of \cite{W7} without proof.
At the end of this section, we discuss the relationships of this
construction with Drinfeld's twistings \cite{Dr6} and twistings of
Landstad and Enock-Vainerman \cite{Landstad1,EnockVain1,Vainerman2}.
In Sect. 4, we construct a non-trivial finite quantum group
$A$ of order $18$ such that neither $A$ nor its dual can be expressed as
a crossed product of the form $C(G_1) \rtimes_{\tau} G_2$, where
$G_1$ and $G_2$ are ordinary finite groups.
Finally, in Sect. 5, we deform finite groups of Lie type
using their maximal abelian subgroups (tori).
\vv
{\em A Convention on Terminology}:
When $A=C(G)$ is a Woronowicz Hopf $C^*$-algebra, we also
call $A$ a compact quantum group, referring to the abstract dual $G$.
Hence a representation of the {\em quantum group} $A$ is a representation
of $G$ in the sense of \cite{Wor5}, which is also called a corepresentation
of the Woronowicz Hopf $C^*$-algebra $A$ (cf. also \cite{W1}); while a
representation of the {\em algebra} $A$ has an obvious different meaning.

\section{
Deformation of algebras via actions of finite
abelian groups}

In this section, we adapt the construction of the monograph of
Rieffel \cite{Rf6} to the situation of actions of finite abelian groups on
$C^*$-algebras (as opposed to actions of ${\Bbb R}^d$ considered there by
Rieffel). Namely, for every quadruple $(A, H, \alpha, J)$ consisting of a
finite dimensional unital $C^*$-algebra $A$, an action $\alpha$ of a
finite abelian group $H$ on the $C^*$-algebra $A$,
and a skew-symmetric automorphism $J$
(with respect to a Pontryagin pairing--see definition below) on $H$,
we construct a deformed unital $C^*$-algebra $A_J$.
It is not our intention to generalize in detail everything in \cite{Rf6}
to this setting.
As a matter of fact, many results in \cite{Rf6} do not generalize
to this setting.
Our primary task in this section is to give some details of those results
that are needed in the next section for the deformation of finite
quantum groups, the main one being the construction of the $C^*$-algebra
$A_J$ mentioned above. We will also briefly indicate some other results that
might be useful elsewhere.

Throughout this section, $A$ will denote a
finite dimensional unital $C^*$-algebra
on which a finite abelian group $H$ acts by $*$-automorphisms $\alpha$.
The group operation of $H$ is written additively. Let
$$H \times H \longrightarrow {\Bbb T}, \; \; \;
(s, t) : \longmapsto <s, t>$$
be a pairing (with values in the circle group ${\Bbb T}$)
that identifies $H$ with its Pontryagin dual $\hat{H}$
(we call such a pairing a {\bf Pontryagin pairing}).
More precisely, identifying $H$ with
$$\ZZ/n_1\ZZ \oplus \ZZ/n_2\ZZ \oplus \cdots
\oplus \ZZ/n_l\ZZ,$$
where $n_1,n_2, \cdots, n_l$ are (not necessarily distinct) natural numbers,
a pairing is given by
\bgeq
<s,t>=<(s_1, \cdots, s_l), (t_1, \cdots, t_l)>
= e^{2 \pi i ( s_1 t_1 /n_1 + s_2 t_2 /n_2 + \cdots + s_l t_l /n_l)},
\ndeq
where $s_k, t_k \in \ZZ, k = 1, \cdots, l$.
%

Let $End(H)$ be the ring of endomorphisms of the group $H$
and $GL(H)$ the group of automorphisms on $H$ (which is
the same as the group of invertible elements in the ring $End(H)$).
Using a Pontryagin pairing $<s, t>$ on $H$ above, we can
define the notion of transpose $J^t$ of an endomorphism $J \in End(H)$.
More generally, if $G$ and $H$ are two finite abelian groups endowed with
Pontryagin pairings, then every group homomorphism $J$ from $G$ to $H$ admits
a transpose $J^t$, which is a homomorphism from $H$ to $G$.
Throughout this section, we assume that $H$ admits a nontrivial
skew-symmetric automorphism.
(Note that some finite abelian groups do not have such automorphisms! But
the examples of groups we consider later in this paper do.)
We can also define the group of orthogonal automorphisms
$O(H)$ in the evident manner.
Just as in the case of a vector space,
by choosing $l$ cyclic generators of $H$, one for
each of the subgroup $H_k \cong \ZZ /n_k\ZZ$ of $H$, we can also represent
elements of $End(H)$ in terms of matrices with entries
consisting of group homomorphisms from $H_j$ to $H_k$, $j, k = 1, \cdots, l$.
With each choice of cyclic generators of $H$,
$GL(H)$ and $O(H)$ 
can be identified with the sets of invertible and orthogonal
matrices, respectively. Note that the
skew-symmetry of the matrix $J$ is independent of the choice of the
generators of $H$.

For any finite group $H$, we will use $\int$ to denote the normalized
Haar integral on $H$, i.e.
\bgeq
\int_{s \in H} f(s) = \frac{1}{|H|} \sum_{s \in H} f(s),
\ndeq
where $|H|$ is the number of elements in $H$ and
$f$ is a function on $H$ taking values in some vector space.
We will see that the normalization is convenient for
the constructions of our deformed algebra and quantum group.
The symbol $\int_{s_1, s_2, ..., s_k \in H}$ will denote the corresponding
$k$-th fold integral.

For the convenience of the reader, we recall here the orthogonality
relations for group characters on $H$, namely the relations
\bgeq
\label{orthogrelation}
\int_{ t \in H} < s, t > = \delta_{s, 0},
\ndeq
where $<s,t>$ is a Pontryagin pairing on $H$.

We will also need the Fourier inversion formula for $A$-valued
functions $F(s)$ on $H$, which we recall also,
as we will be using it a number of times,
\bgeq
\label{FourierInversion1}
\Check{\Hat{F}} = F, \; \; i.e., \; \;
|H| \int_{s, t \in H} F(s) <s, -t> <t, x> = F(x).
\ndeq
(The inversion formula is easily seen to be a consequence of
the orthogonality relations \eqref{orthogrelation}.) In
particular, we have,
\bgeq
\label{FourierInversion2}
|H| \int_{s, t \in H} F(s) <s , t> = F(0).
\ndeq

\bgdf
\label{deformedproduct}
(cf \cite{Neumann1,Rf6})
Let $J \in End(H)$ be an endomorphism on $H$.
The {\bf deformed product} $\times_{J}$ (or $\times_J^\alpha$) on $A$ is
defined by
\bgeq
a \times_{J} b = |H| \int_{s,t \in H}
\alpha_{s}(a) \alpha_{t}(b) <J s,t>, \; \; \; a,b \in A,
\ndeq
where the products on the right hand side is in $A$.
Let $A_J$ (or $A_J^\alpha$)  denote $(A, \times_J)$.
\nddf

The number $|H|$ in the above formula insures that the deformed
algebra $A_J$ is unital (see (2) of the next proposition).
The observant reader might have noticed that in the above formula
we have chosen $J$ to appear in the dual pairing $< \cdot, \cdot >$
instead of in the action $\alpha$, as is done in \cite{Rf6,Rf8,W4}
(cf also von Neumann \cite{Neumann1}).
We do this because if we replace $J$ with $\JJ$, where $\hbar$ is any
{\em real} number, the above formula still make sense. But $\alpha_{\JJ s}$
does not make sense for finite group $H$ acting on $A$. See also
the remarks at the end of this section.

\bgprop
\label{associativity}
(1). For any $J \in End(H)$, the deformed product $\times_{J}$ is associative.

(2). If $J \in GL(H)$, then the unit of $A$ continues to be the unit of $A_J$.
\ndprop
\pf
(1).
This is the analog of Theorem 2.14 of Rieffel \cite{Rf6}. However the proof
given there does not work for finite abelian groups because their subgroups
do not have orthogonal complements. We give a much simpler proof of this
result. The key is to make the {\em correct} change of variables.
We compute, by applying change of variables twice
(the second change of variable is a little bit tricky),
\bgeqq
( a \times_{J} b ) \times_{J} c &=&
|H| \int_{s,t \in H}
\alpha_{s}(a \times_{J} b) \alpha_{t} (c) <J s,t> \\
&=& |H|^2 \int_{s,t, u, v \in H}
\alpha_{s + u} (a) \alpha_{s + v} (b) \alpha_{t} (c) <J s,t> <J u, v> \\
&=& |H|^2 \int_{s,t, u, v \in H}
\alpha_{s} (a) \alpha_{s - u + v} (b) \alpha_{t} (c)
<J (s - u), t> <J u, v> \\
&=& |H|^2 \int_{s,t, u, v \in H}
\alpha_{s} (a) \alpha_{u + v} (b) \alpha_{t}(c) <J u , t> <J (s - u), v>,
\ndeqq
which, after exchanging the roles of $t$ and $v$,
\bgeqq
&=& |H|^2 \int_{s,t, u, v \in H}
\alpha_{s} (a) \alpha_{u + t} (b) \alpha_{v} (c) <J u , v>
<J (s - u), t> \\
&=& |H|^2 \int_{s,t, u, v \in H}
\alpha_{s} (a) \alpha_{u + t} (b) \alpha_{v} (c)
<J u , v - t> <J s , t> \\
&=& |H|^2 \int_{s,t, u, v \in H}
\alpha_{s} (a) \alpha_{u + t} (b) \alpha_{v + t} (c) <J u , v> <J s , t>.
\ndeqq
On the other hand, expanding
$ a \times_J (b  \times_J c) $ we see that
$$ a \times_J (b  \times_J c) =
|H|^2 \int_{s,t, u, v \in H}
\alpha_{s} (a) \alpha_{t + u} (b) \alpha_{t + v} (c) <J u , v> <J s , t>.$$

(2). If $J \in GL(H)$, then $<J s, t>$ is a Pontryagin pairing for $H$, hence
we can use \eqref{FourierInversion2} (replacing $<s, t>$
in \eqref{FourierInversion2} by $<J s, t>$),
\bgeqq
a \times_{J} 1
= |H| \int_{s,t \in H} \alpha_{s}(a) <J s,t>
= \alpha_{0} (a) = a.
\ndeqq
Similarly,
$$1 \times_J b = b.$$
This proves the proposition.
\QED
\vv
As in \cite{Rf6}, let $A_u$ be the spectral subspace of $u \in H$:
\bgeq
A_u = \{ a \in A \; | \; \alpha_s(a) = <u , s> a, s \in H \}.
\ndeq

\bgprop
\label{specsubspace}
Let $J \in GL(H)$ be a skew-symmetric automorphism: $J^t = -J$.
Let $a \in A_u$, $b \in A_v$ (the spectral subspace of $A$
corresponding to $u, v \in H$). Then
\bgeq
a \times_J b = <J^{-1} u, v> a b.
\ndeq
\ndprop
\pf
This is the analog of Proposition 2.22 in \cite{Rf6}.
Instead of the Poisson summation formula, as is used in
the proof of 2.22 in \cite{Rf6}, we apply the Fourier
inversion formula to the last line of the following computation:
\bgeqq
a \times_{J} b
&=& |H| \int_{s,t \in H} \alpha_{s}(a) \alpha_{t}(b) <J s,t> =
    |H| \int_{s,t \in H} <s, u> a <t, v> b <J s,t>
    \\
&=&
|H| \int_{s,t \in H} <J^{-1} s, u> <t, -v> <s, -t> ab = <J^{-1}(-v), u> a b.
\ndeqq
That is
$$a \times_{J} b = <J^{-1} u, v> a b$$
by skew-symmetry of $J$.
\QED
\vv
{\em Remark.}
Note that if $A$ is commutative, then
$a \times_J b = <2J^{-1} u, v> b \times_J a,$
where $a, b$ are as in the above proposition.
Hence, we see that if $2 J^{-1} \neq 0$ and if the action $\alpha$ is
non-trivial, then the algebra $A_J$ is noncommutative,
even if $A$ is a commutative algebra.
The condition $2J^{-1} = 0$ is related with the characteristic
$2$ phenomenon (see the last two sections for examples
concerning this).

\bgprop
Let $J \in End(H)$ be a skew-symmetric homomorphism: $J^t = -J$. Then
under the involution $*$ of the algebra $A$, we have
$$(a \times_J b)^* = b^* \times_J a^*$$
for $a, b \in A_J$. Hence $A_J$ is a $*$-algebra.
\ndprop
\pf
Use $<Jy , x> = <y , J^t x> = <y , -Jx>$.
\QED
\vv
Consider the Hilbert $A$-module $E = C(H) \otimes A$
under the $A$-valued inner product
\bgeq
<f, g>_A = \int_{x \in H} f^*(x) g(x), \; \; f, g \in C(H) \otimes A,
\ndeq
where $C(H)$ is the algebra of complex valued functions on $H$.
Note that as a tensor product of two $C^*$-algebras, $E$ is also a
$C^*$-algebra and $H$ acts on it by translation:
\bgeq
\tau_s(f)(x) = f(x-s).
\ndeq
If $J$ is a skew-symmetric automorphism, then from the propositions
above, $E_J = (E, \times^\tau_J)$ is a unital $*$-algebra. Let $L$ denote the
left regular multiplication on $E_J$:
\bgeq
L_f g = f \times^\tau_J g.
\ndeq
Under these assumptions, we have (cf 4.2, 4.3, 4.6 of \cite{Rf6})

\bgprop
The left regular multiplication $L$ is a faithful unital $*$-representation of
the $*$-algebra $E_J$ by bounded operators on the Hilbert $A$-module $E$. More
precisely, we have
$ L_f = 0$ if and only if $f=0$,  and
\bgeq
\label{L*rep}
<f \times_J^\tau g, h>_A = < g, f^* \times_J^\tau
h>, \; \; f \in E_J, \; g, h \in E,
\ndeq
 \bgeq
\label{Lnorm}
||L_f|| \leq \int_{s \in H} ||f(s)|| = ||f||_1, \; \; \; f \in E_J.
\ndeq
\ndprop
\pf
The identity \eqref{L*rep} is a straightforward checking without going into
the complications such as involved in the proof of 4.2 in
Rieffel \cite{Rf6}. We leave this to the reader.

We show that $L$ is faithful.
Let $L_f=0$. Hence
$$<f \times_J^\tau g, f \times_J^\tau g>_A =0, \; \; \; g \in E.$$
Let $g$ be the unit element of the $C^*$-algebra $E$, $g(s) = 1, \; \; s
\in H$. Then by \propref{associativity}, $g$ is the unit of $E_J$. Hence
$$<f \times_J^\tau g, f \times_J^\tau g>_A = <f, f>_A = 0.$$
Note that $g$ plays two roles in here: as the unit of the {\em algebra} $E_J$
and as an vector in the Hilbert $A$-{\em module} $E$.
Hence $f = 0$.

The proof of the inequality \eqref{Lnorm} is the same as (and easier than)
the proof of 4.3 in \cite{Rf6} (see also 4.6 of \cite{Rf6}).
For the convenience of the reader, we sketch the proof here.
A short computation shows that
$$
L_f(g) = f \times_J^\tau g = \int_{s \in H} f(s) U_s(g),
\; \; i.e., \; \;  L_f = \int_{s \in H} f(s) U_s,
$$
where $U_s$ is the unitary operator on the Hilbert module $E$ defined by
$$U_s(g)(x) = <J(x-s), x> \Check{g}(J(x-s)),$$
$\Check{g}$ being the {\em inverse} Fourier transform of $g$
(Plancherel's theorem). Hence \eqref{Lnorm} is immediate.
\QED
\vv
Let us come back to our algebra $A_J$. For $a \in A_J=A$, the element
$\tilde{a}$ of $E_J$ defined by $\tilde{a}(s) = \alpha_s(a)$
is zero if and only if $a=0$. Using the above result, we can define
a $C^*$-norm on $A_J$ as follows.
\bgdf
\label{deformednorm}
Let $J$ be a skew-symmetric automorphism on $(H, <,>)$.
The deformed $C^*$-norm $||\cdot||_J$ on $A_J$ is defined by
\bgeq
||a||_J = ||L_{\tilde{a}}||, \; \; a \in A_J,
\ndeq
where $||L_{\tilde{a}}||$ is the operator norm of $L_{\tilde{a}}$ on the
Hilbert $A$-module $E$.
\nddf

Summarizing the above, we have the following main result of the section:
\bgth
\label{theoremalgebra}
Let $A$ be a finite dimensional
unital $C^*$-algebra. Let $H$ be a finite abelian group
acting on $A$ by automorphisms. Let $J \in GL(H)$ be a skew-symmetric
automorphism: $J^t = -J$. Then $A_J$ is a unital $C^*$-algebra under the norm
$||\cdot||_J$.
\ndth
\noindent
{\em Remarks.}
(1). Note that on any finite dimensional $*$-algebra,
there can be at most one $C^*$-norm. Hence the $C^*$-norm
defined above is the unique one on $A_J$.

(2). Note that we need $J$ to be a skew-symmetric automorphism in order to
define the $C^*$-norm on $A_J$, while in Rieffel \cite{Rf6}, $J$
can be any skew-symmetric endomorphism on a vector space.
Also note that we do not have the analog of
Theorems 2.15 and 6.5 of Rieffel \cite{Rf6}. Namely,
$$(A_J)_K = A_{J + K}$$
is not true in general. However, using the orthogonality
relations for group characters (see \eqref{orthogrelation}),
one can easily prove the following proposition.

\bgprop
\label{deformback}
Under the assumption of the theorem above, $(A_J)_{-J} = A$.
\ndprop

Now we can state the analogs of Theorems 2.10. 5.7, 5.8, 5.12 and 7.7 of
Rieffel \cite{Rf6}.

\bgprop
\label{exactsequence}
Let $J \in GL(H)$ be a skew-symmetric automorphism.
Let $\alpha$ and $\beta$ be
actions of $H$ on $A$ and $B$ respectively.
Let $\theta: A \longrightarrow B$ be an equivariant homomorphism.

(1). $\theta$ is still an equivariant homomorphism from $A_J$ to $B_J$
(denote this homomorphism by $\theta_J$, called the deformation of $\theta$);

(2). $\theta$ is injective (resp. surjective) if and only if $\theta_J$ is.

(3). Let $I$ be an ideal of $A$ that is invariant under the action $\alpha$.
Let $Q = A/I$, and let $\alpha$ also denote the action of $H$ on $Q$, so
we have an equivariant exact sequence
$$0 \longrightarrow I \longrightarrow A \longrightarrow Q \longrightarrow 0.$$
Then the corresponding sequence (see (1) above)
$$0 \longrightarrow I_J  \longrightarrow A_J
\longrightarrow Q_J \longrightarrow 0$$
is also exact.
\ndprop
The proofs of these analogs follows directly from our {\em definitions}
and the key assumption that $A$ is a finite dimensional $C^*$-algebra.
We leave the checking to the reader.
The reader is advised not to look up the proofs in \cite{Rf6} for clues
(for otherwise the reader would be mislead to complicate the proofs of
these analogs), but to simply think about our definitions.
\vv
{\em Remarks.}
(1). If we replace $J$ with $\JJ$ in the construction above,
a number of things in this section are still true, where $\hbar$ is real, and
 $$<\hbar J s, t> =
= e^{2 \pi \hbar i ( s_1 t_1/n_1 + s_2 t_2/n_2 + \cdots + s_l t_l/n_l)},$$
using the above identification of $H$ with the concrete abelian group
as a direct sum of cyclic groups. For any skew-symmetric $J$,
$A_\JJ$ is an associative $*$-algebra, but it may not have a unit
or a $C^*$-norm even if $J$ is an automorphism.
So it not clear how one constructs strict deformation quantization.

(2). For practical purposes of next section, we have restricted $A$ to be
a finite dimensional $C^*$-algebra. If we remove this restriction,
then the proofs of all the above results, except \propref{exactsequence},
are still valid (Of course, in \thref{theoremalgebra}, we need
a completion to obtain a $C^*$-algebra).
We believe that \propref{exactsequence} is still true in this case.

\section{
Deformation of finite quantum groups via finite
abelian subgroups}

In the theory of finite groups,
the finite groups of Lie type 
are one of the most important classes of finite groups.
In view of the fact that classical Lie groups have $q$-deformation, a
natural question in this connection is
\bgprob
\label{finitegroupsoflietype}
Do finite groups of Lie type have
an analog of $q$-deformation into finite quantum groups?
\ndprob

This problem seems to be out of reach at the moment.
In this section, we construct a deformation of Rieffel type for finite groups
(as well as for finite quantum groups) that contain an abelian subgroup.
This deformation is not the analog of the $q$-deformation,
it is dual to Drinfeld twistings of the
quantized universal enveloping algebras.
We will see this at the end of this section.

We start with a finite quantum group $G =(A, \Phi)$ (in the sense
that $A$ is the ``function space'' $C(G)$, where $\Phi$
is the coproduct on $A$ \cite{Wor5}). Assume that its maximal subgroup
$X(A)= \{$ *-homomorphisms from $A$ into ${\Bbb C} \}$
contains an abelian subgroup $T$ with a
nontrivial skew-symmetric automorphisms $S$.
So there is a surjective morphism of Hopf $C^*$-algebras
$\pi: A \longrightarrow C(T)$.
Let
\bgeq
H: = T \oplus T,
\ndeq
and let
\bgeq
J: = S \oplus (-S)
\ndeq
be the skew-symmetric automorphism on $H$.
Define an action $\alpha$ of $H$ on the $C^*$-algebra $A$ as follows:
\bgeq
\label{TheAction}
\alpha_{(s,u)} = \lambda_{s} \rho_{u},
\ndeq
where
\bgeq
\lambda_{s} = (E_{-s} \pi \otimes id) \Phi, \; \; \;
  \rho_{u} = (id \otimes E_{u} \pi ) \Phi,
\ndeq
$id$ being the identity map on $A$ and $E_u$ the evaluation functional
on $C(T)$ corresponding to $u$.
Using results of the previous section, we obtain a deformed $C^*$-algebra
$A_J$ with new product $\times_J$ defined by
(see  formula \eqref{deformedproduct})
\bgeq
\label{deformedproduct1}
a \times_J b = |T|^2 \int_{s,t,u,v \in T}
\alpha_{(s,u)}(a) \alpha_{(t,v)}(b)
<Ss, t><-S u, v>,
\ndeq
where $a,b \in A$ and  $<s,t>$ is a Pontryagin pairing on $T$.
The main result of this section is (cf \cite{Rf8,W4})
\bgth
\label{theoremquantumgroup}
Under the same coproduct $\Phi$ of
$A$, the deformation $(A, \times_J)$ is still a
finite quantum group containing $T$ as a subgroup.
\ndth
\noindent
{\em Remarks on the proof.}
We will show that $A_J$ satisfies the axioms of a finite dimensional
Hopf $C^*$-algebra as given in Kac and Palyutkin \cite{KP2},
instead of the ones given in Appendix 2 of Woronowicz
\cite{Wor5}, though they are equivalent to each other.
The proof is a modification of the proof of Theorem 3.9 of \cite{W4}.
Unlike that theorem, because $A$ is of finite dimension here,
we do not need to consider
the analogs of Propositions 3.2 and 3.8 in \cite{W4}, which are
essential steps for the proof of that theorem.
The subtlety in our situation here is that
the method used there in the treatment of the deformed coproduct does not
work anymore, because the existence of
orthogonal complements is used in an essential way there (but,
as pointed out before, a subgroup of a finite abelian group
needs not have an orthogonal complement).
To deal with the deformed coproduct, we will show that the heuristic
computation on page 471 of \cite{Rf8} can be made rigorous in
our setting (replacing the compact Lie group $G$ there by
our finite quantum group).
\pf
Let $F, G \in A_J \otimes A_J$, and let $\times_J$ also denote the product
on $A_J \otimes A_J$.
Using formula \eqref{deformedproduct} we can find the formula
for the product in the $C^*$-algebra $A_J \otimes A_J$ in terms of product in
the $C^*$-algebra $A \otimes A$, with
the summation (integration)  over repeated indices:
\bgeq
\label{deformedtensor}
F \times_J G  &=& |T|^4 \int 
                 \gamma_{(s,u, s',u')}(F)
                 \gamma_{(t, v, t', v')}(G)
                  <L(s, u, s', u'),(t, v, t', v')>
\ndeq
where $\gamma = \alpha \otimes \alpha$ is the tensor product action of
$H \oplus H$ on $A \otimes A$, $L = J \oplus J$ is the corresponding
skew-symmetric automorphism on $H \oplus H$,  and
\bgeq
<(s, u, s', u'),(t, v, t', v')> = <s, t> <u, v> <s', t'> <u', v'>.
\ndeq
Note that this identity is easily verified on tensors of the form
$$F= a_1 \otimes a_2, \; \; \; G= b_1 \otimes b_2.$$
Since $A$ is finite dimensional, this gives an isomorphism of
$*$-algebras
$$(A \otimes A)_L^\gamma = A_J^\alpha \otimes A_J^\alpha.$$
It is easy to see this is actually an isomorphism of $C^*$-algebras
(see Remark (1) after \thmref{theoremalgebra}).
This isomorphism is the analog of Corollary 2.2 of \cite{Rf8}, where a
more complicated proof is needed.

For $a, b \in A_J$, we have by \eqref{deformedtensor} and \eqref{TheAction}
(cf \cite{Rf8})
\bgeqq
\Phi(a) \times_J \Phi(b) &=&
\Phi(a) \times_L \Phi(b) \\
&=& |T|^4 \int_{s,u,t,v,s',u',t',v' \in T}
    (\lambda_{s} \rho_{u} \otimes \lambda_{s'} \rho_{u'})(\Phi(a))
    (\lambda_{t} \rho_{v} \otimes \lambda_{t'} \rho_{v'})(\Phi(b))
    \\
& & \hspace{2cm}
    <Ss, t> <-Su, v> <Ss', t'> <-Su', v'>,
\ndeqq
which, by 2.7 of \cite{W4}
\bgeqq
&=& |T|^4 \int_{s,u,t,v,s',u',t',v' \in T}
    (\lambda_{s} \rho_{u-s'} \otimes  \rho_{u'})(\Phi(a))
    (\lambda_{t} \otimes \lambda_{t'-v} \rho_{v'})(\Phi(b))
    \\
& & \hspace{2cm}
    <Ss, t> <-Su, v> <Ss', t'> <-Su', v'>.
\ndeqq
Making change of variables $u-s' \mapsto u$, $t'-v \mapsto t'$, and using
\eqref{FourierInversion2} twice (note that both $<-Su, v>$ and $<Ss', t'>$ are
Pontryagin pairings on $T$!), the last expression
\bgeqq
&=& |T|^4  \int_{s,u,t,v,s',u',t',v' \in T}
    (\lambda_{s} \rho_{u} \otimes  \rho_{u'})(\Phi(a))
    (\lambda_{t} \otimes \lambda_{t'} \rho_{v'})(\Phi(b))
    \\
& & \hspace{2cm}
    <Ss, t> <-Su, v> <Ss', t'> <-Su', v'>
    \\
&=& |T|^2  \int_{s,t,u',v' \in T}
    (\lambda_{s} \otimes  \rho_{u'})(\Phi(a))
    (\lambda_{t} \otimes  \rho_{v'})(\Phi(b))
    <Ss, t> <-Su', v'>,
\ndeqq
which, by 2.7 of \cite{W4}
\bgeqq
&=& |T|^2  \int_{s,t,u,v \in T}
   \Phi( \lambda_{s}  \rho_{u} (a)) \Phi( \lambda_{t}  \rho_{v} (b))
    <Ss, t> <-Su, v>
    \\
&=& \Phi(a \times_J b).
\ndeqq
That is
$$\Phi(a) \times_J \Phi(b) = \Phi(a \times_J b).$$

As in \cite{Rf8}, the action $\alpha$ restrict to an action
on $C(T)$ and $\pi$ is equivariant. From \propref{exactsequence} of
the last section, this gives
a surjective homomorphism $\pi_J$ from $A_J$ onto $C(T)_J$.
It is also clear that
\bgeq
(\pi_J \otimes \pi_J) \Phi_J = \Phi_T \pi_J,
\ndeq
where $\Phi_T$ is the coproduct on $C(T)$.
However the method used in \cite{Rf8} for the proof of
\bgeq
C(T)_J = C(T)
\ndeq
does not work here, because it uses a result of \cite{Rf6}
which is not true for finite abelian groups. We can prove this
directly as follows. For $f \in C(T)$, we have
$$\alpha_{(s,u)}(f) = \lambda_s \rho_u (f) = \lambda_{s-u}(f).$$
Hence
\bgeqq
f \times_J g  &=& |T|^2  \int_{s,u,t,v \in T}
                  \lambda_{s-u}(f) \lambda_{t-v}(g) <Ss, t> <-Su, v>
                  \\
              &=& |T|^2  \int_{s,u,t,v \in T}
                  \lambda_{s}(f) \lambda_{t}(g) <Ss, t> <Ss, v> <Su, t>  \\
               &=& f g,
\ndeqq
where we have used the orthogonality relations of the characters of
a finite abelian group.
This shows that $T$ will still be a subgroup of $A_J$ once $A_J$ is
shown to be a quantum group.

The counit of $A_J$ is defined by
\bgeq
\epsilon_J = \epsilon_T \pi_J,
\ndeq
where $\epsilon_T$ is the counit of $C(T)$. So as a linear map,
$\epsilon_J$ is the same as $\epsilon$.
The coinverse $\kappa_J$ on $A_J$ is defined to be the same as $\kappa$.
The identity
\bgeq
\kappa_J (a \times_J b) = \kappa_J(b) \times_J \kappa_J(a)
\ndeq
is a direct consequence of the fact that
$$\kappa \alpha_{(s,u)} = \alpha_{(u,s)} \kappa$$
(cf 2.8 of \cite{W4}) and the skew-symmetry of $S$.

Now we check the antipodal property
\bgeq
m_J(id_J \otimes \kappa_J) \Phi_J
= I_J \epsilon_J
=m_J(\kappa_J \otimes id_J) \Phi_J ,
\ndeq

By 2.8 and 2.6  of \cite{W4}, we have for coefficients
$a_{ij} \in A$ of a unitary representation $(a_{ij})$
of the quantum group $A$ that
\bgeqq
& & m_J(id_J \otimes \kappa_J) \Phi_J (a_{ij}) =
    m_J(id_J \otimes \kappa_J) \Phi (a_{ij})
    \\
&=& |T|^2
\int_{s,u,t,v \in T} m(\alpha_{(s,u)} \otimes \alpha_{(t,v)} \kappa )
\Phi (a_{ij}) <J(s,u), (t,v)> \\
&=& |T|^2   \int_{s,u,t,v \in T}  m(id \otimes \kappa)
(\lambda_{s} \rho_{u} \otimes \lambda_{v} \rho_{t})
\Phi(a_{ij}) <Ss, t><-Su, v>  \\
&=& |T|^2   \int_{s,u,t,v \in T}  m(id \otimes \kappa)
( \lambda_s \rho_{u-v} \otimes \rho_{t})
\Phi(a_{ij}) <Ss, t><-Su, v>  \\
&=& |T|^2   \int_{s,u,t,v \in T}  m(id \otimes \kappa)
( \lambda_s \rho_{u} \otimes \rho_{t})
\Phi(a_{ij}) <Ss, t><-Su, v>,
\ndeqq
which, using \eqref{FourierInversion2}, noting that
$<-Su, v>$ is a Pontryagin pairing on $T$,
\bgeqq
&=& |T|   \int_{s,t \in T}  m(id \otimes \kappa)
( \lambda_s  \otimes \rho_{t})
\Phi(a_{ij}) <Ss, t> \\
&=& |T|   \int_{s,t \in T}  m ( \lambda_s  \otimes \kappa \rho_{t})
(\sum_k a_{ik} \otimes a_{kj}) <Ss, t> \\
&=& |T|   \int_{s,t \in T}
\sum_{k,l,r} E_{-s}(\pi(a_{il})) a_{lk} a^*_{rk} E_t(\pi(a_{rj})) <Ss, t> \\
&=& |T|   \int_{s,t \in T}
(E_{-s} \otimes E_t) \Phi_T (\pi(a_{ij})) <Ss, t> \\
&=& |T|   \int_{s,t \in T}
E_{-s + t}  (\pi(a_{ij})) <Ss, t>
= |T|   \int_{s,t \in T}
E_{t}  (\pi(a_{ij})) <Ss, t> ,
\ndeqq
which, using once again \eqref{FourierInversion2},
\bgeqq
&=& E_0 (\pi(a_{ij})) = \epsilon_T(\pi(a_{ij}))=
\epsilon(a_{ij}) = \epsilon_J (a_{ij}).
\ndeqq
That is, on $A_J$ (note that $A_J = A$ as a vector space),
$m_J(id_J \otimes \kappa_J) \Phi_J = I_J \epsilon_J.$

Similarly $m_J(\kappa_J \otimes id_J) \Phi_J = I_J \epsilon_J$.
\QED
\vv
We will denote the coproduct of $A_J$ by $\Phi_J$.
Note that because of \propref{deformback},
$A_J$ can be deformed back to the original
quantum group $(A_J)_{-J} = A$.
\vv
{\em Remarks.}
(1). In the above,
we assumed that $A$ is a finite quantum group instead of a more general
compact quantum group. In the latter case, we do not know how to show
directly that $A_J$ is still a compact quantum group. The main difficulty
is to rigorously define the coproduct. Note that
we can not define an analog of the
map $ \varrho $ as given near 3.3 of \cite{W4} (see also \cite{Rf8}).
One possible approach is to use the Krein duality as given in our paper
\cite{W3}. 

(2). Just as in \cite{Rf8,W4}, the Haar measure of $A$ is still
the Haar measure of $A_J$. One can easily see this from the
uniqueness of the (left and right invariant) Haar measure and the fact that
the coproduct of $A_J$ is that of $A$. This proof is
also valid for the infinite dimensional cases of \cite{Rf8,W4} if we
work with the dense Hopf $*$-algebra, and it is much easier
than the proof in \cite{Rf8}.

(3). From the above remark and the orthogonality relations for
characters of irreducible representations,
we see that the irreducible representations
of the quantum group $A$ are still irreducible representations of
the quantum group $A_J$.
From these, we see that the representation ring of the quantum group $A$
is invariant under deformation, just as in \cite{Rf8,W4} and \cite{Nik}.
\vv
Now we describe the construction above in the dual picture.
Let $B$ be a finite dimensional Hopf $C^*$-algebra with coproduct $\Phi$.
Let $T$ be an abelian subgroup of the group of group-like elements of $B$.
Endow $T$ with a Pontryagin pairing. For notational
convenience, we now assume that the group operation on $T$
is multiplicative instead of additive, and for a skew
symmetric automorphism $S$ on $T$, the matrix $-S$ will denote the
automorphism $(-S) x = (Sx)^{-1}$.
As before, let $J = (S, -S)$ be the skew-symmetric automorphism
on $H = T \times T$. Then $B$ is a Hopf $C^*$-algebra
under the original $C^*$-algebra structure, original counit
and antipode, and the deformed coproduct given by
\bgeq
\Phi_J (b) = |T|^2 \int_{s,u,t,v \in T}  (\beta_{(s,u)} \otimes \beta_{(t,v)})
\Phi (b) <Ss, t> <(Su)^{-1}, v>,
\ndeq
where
\bgeq
\beta_{(s,u)}(b) = s b u^{-1}.
\ndeq
Using orthogonality relations for group characters,
we have that
\bgeq
(|T| \int_{s,t \in T}  (s \otimes t) <Ss, t>)^{-1} =
|T| \int_{u,v \in T} (u^{-1} \otimes v) <Su, v>.
\ndeq
We can also check that
\bgeq
\label{pseudo-2-cocycle}
F = |T| \int_{s,t \in T} (s \otimes t) <Ss, t>
\ndeq
is a unitary element in $B \otimes B$.
We can now rewrite $\Phi_J(b)$ as
\bgeq
\Phi_J (b) =  F \Phi (b) F^{-1}.
\ndeq
This shows that our deformation of finite quantum groups
in this dual picture is an analog of the twistings of quantized
universal enveloping algebra of Drinfeld \cite{Dr6,Lev1,LS1}, just as
the ones in \cite{Rf8,W4}. The element $F$ is not a 2-cocycle,
though both $(\epsilon \otimes 1) (F) =1$ and $(1 \otimes \epsilon) (F) =1$
are satisfied. But since the twist $\Phi_J$ is coassociative because of
\thref{theoremquantumgroup}, $F$ is a pseudo-2-cocycle in the sense of
\cite{EnockVain1,Vainerman2}.
However, it is not clear how to verify directly that $F$
satisfies the pseudo-2-cocycle condition for a general noncommutative and
noncocommutative Hopf $C^*$-algebra $B$ without passing to the dual
$A$ of $B$ (see \propref{associativity} and \thref{theoremquantumgroup})!
The above twisting of $B$ can be viewed as
the {\em canonical} one among the ones considered in
\cite{EnockVain1,Vainerman2} that are associated with
a finite abelian subgroup $T$.
We will call the twist $F$ the {\em Weyl-von Neumann-Rieffel twist}
associated with $(T, S)$, and denote the twist of $B$ by $B^S$ (to distinguish
it from $A_J$). To gain a better understanding of this construction,
it is desirable to solve the following (cf remark (3) above and \cite{Nik})
\bgprob
Characterize the finite quantum groups that can be
obtained as deformations of finite groups in the manner above.
Find their isomorphic invariants.
\ndprob

\section{
A finite quantum group of order $18$}

Let $T = ({\Bbb Z}/n{\Bbb Z})^{2k} $ with canonical
generators $\epsilon_j$, $j = 1, \cdots, 2k$ ($\epsilon_j$ is the
element with $j$-th component $1$ and the other components $0$).
Let $\beta$ be the action of ${\Bbb Z}/2{\Bbb Z}$ on $T$ which
exchanges $\epsilon_j$  and $\epsilon_{j + k}$, $j = 1, \cdots, k$.
Let $G$ be the corresponding semi-direct product of $T$ by
${\Bbb Z} / 2{\Bbb Z}$ under this action:
$G = T \rtimes_{\beta} {\Bbb Z} / 2 {\Bbb Z}$. Let $A = C(G)$.
Consider the skew-symmetric automorphism $S$ of $T$ defined by
\bgeq
S(\epsilon_j) = - \epsilon_{k + j}, \; \; \;
S(\epsilon_{k + j}) =  \epsilon_{j}, \; \; \;  j = 1, \cdots, k.
\ndeq
Alternatively, $S$ has matrix representation
$$  \left(
\begin{array}{cccc}
0 & I_k    \\
-I_k & 0
\end{array}
\right)
$$
with respect to the
generators (or ``basis'') $\epsilon_j$, $j = 1, \cdots, 2k$,
where $I_k$ represents the identity transformation
on the group $({\Bbb Z}/n{\Bbb Z})^{k} $. Then we are in the position to
apply \thmref{theoremquantumgroup} to obtain a noncommutative
deformation $A_J$, where $J = S \oplus (-S)$
(see the remark after \propref{specsubspace}).

Take $k = 1$ and $n = 3$. Then $A_J = C(G)_J$ is a finite
quantum group of order $18$, where by definition,
the order of a finite group is the dimension of its function algebra.
From the remark after \propref{specsubspace},
we see that this quantum group is noncommutative and noncocommutative.

We now show that the quantum group $A_J$ (and its dual $B^S$) is not
a crossed product of the form $C(G_1) \rtimes_{\tau} G_2$, where $G_1$
and $G_2$ are ordinary finite groups and $\tau$ is an
action of $G_2$ on $C(G_1)$ by automorphisms of quantum groups
(the triple $(C(G_1), G_2 , \tau)$ is
also called a Woronowicz Hopf $C^*$-dynamical system, see \cite{KP2,W2}).

An easy application of the Mackey Machine shows that
the group $G= ({\Bbb Z}/3{\Bbb Z})^{2} \rtimes_{\beta} {\Bbb Z}/2 {\Bbb Z}$
has 6 irreducible representations: 2 one-dimensional representations
and 4 two-dimensional representations.
Hence by remark (3) after the proof of \thref{theoremquantumgroup},
the quantum group $A_J$ also has 6 irreducible representations.
If $A_J$ is a crossed product of the form $C(G_1) \rtimes_\tau G_2$,
then we only need to consider four cases:
(i). $|G_1| = 2$, $|G_2| = 9$,
(ii). $|G_1| = 9$, $|G_2| = 2$,
(iii). $|G_1| = 3$, $|G_2| = 6$,
(iv). $|G_1| = 6$, $|G_2| = 3$.
We claim that in each of the cases (i), (ii) and (iii), the quantum group
$A_J = C(G_1) \rtimes_\tau G_2$
is a group $C^*$-algebra of an ordinary group and hence it
has 18 irreducible representations instead of 6
(cf. \cite{Wor5}, the irreducible
representations of a compact quantum group of the form $C^*(\Gamma)$
are exactly the elements of of the discrete group $\Gamma$).
This claim implies that cases (i), (ii) and (iii) cannot happen.
To prove the claim we note
first that $G_1$ is abelian in each of these three cases. Since
$\tau$ is assumed to preserve the Hopf $C^*$-algebra structure
of $C(G_1)$, which is isomorphic to $C^*(\hat{G}_1)$
(as a Hopf $C^*$-algebra), by transport of structure, $\tau$ is
an action of $G_2$ by automorphisms on the group $\hat{G}_1$. Hence
$C(G_1) \rtimes_\tau G_2 \cong C^*(\hat{G}_1 \rtimes_\tau G_2) $ as claimed.
Now consider case (iv). If $G_1$ is abelian, then
the same reasoning as above leads to a contradiction.
If $G_1$ is non-abelian, then $G_1$ is the group $S_3$, hence it
has 3 irreducible representations. From the classification
of irreducible representations of quantum groups
associated with crossed products (see Theorem 3.7 of \cite{W2}),
$C(G_1) \rtimes_\tau G_2$ has $3 \times 3 = 9$ irreducible representations.
This again contradicts the fact that the quantum group $A_J$ has
$6$ irreducible representations. This shows that case (iv) cannot happen.

Similarly, we show that the dual $B^S$ of the above $A_J$ cannot be
expressed in the form $C(G_1) \rtimes_{\tau} G_2$ either.
As above we only need to look at case (iv) above with $G_1 = S_3$ and $\tau$
non-trivial. By remark (3) after the proof of \thref{theoremquantumgroup}
and the analysis in the previous paragraph,
the algebra $B^S$ has 6 irreducible representations (which are exactly
the irreducible representations of the quantum group $A_J$).
By transport of structure, $\tau$ gives an automorphism
of order three of the group $S_3$.
An examination of the structure of $S_3$ then shows that the
action $\tau$ is conjugation by a 3 cycle in $S_3$.
From this, using Mackey Machine we see that the
algebra $C(S_3) \rtimes_{\tau} {\Bbb Z}/3 {\Bbb Z}$ has 10 irreducible
representations instead of 6: 9 one-dimensional representations and
1 three-dimensional representation. This completes the proof.
\vv
{\em Remarks.}
(1). The group
$D_4 = ({\Bbb Z}/2{\Bbb Z})^{2} \rtimes_{\beta} {\Bbb Z}/2 {\Bbb Z}$
is the only non-abelian group of order $8$ that has a maximal abelian subgroup
$T$ isomorphic to $({\Bbb Z}/2{\Bbb Z})^{2}$, and $S$ defined above is the
only non-trivial skew-symmetric automorphism on $T$. It is
easy to see that $C(D_4)_J = C(D_4)$
(see the remark after \propref{specsubspace}).
Since the $8$ dimensional quantum group of Kac-Palyutkin \cite{KP2}
clearly contains $T$ as a subgroup, we conclude from
\thref{theoremquantumgroup} that this quantum group
is not of the form $C(G)_J$ for any finite group $G$.

(2). By the same method as above, we note that the duals of the two examples
of 12 dimensional Kac algebras in \cite{EnockVain1,Nik}
(as twists of $D_6$ and $Q_{4}$ respectively) are not isomorphic to
a Kac algebra of the form $C(S_3) \rtimes_{\tau} {\Bbb Z}/ 2{\Bbb Z}$
(like case (iv) above, this is the only non-trivial case),
where ${\Bbb Z}/ 2{\Bbb Z}$ acts on $S_3$ non-trivially.
The second tensor power of one of the 2 two-dimensional irreducible
representations of the algebra $L(G)$ of \cite{EnockVain1}
decomposes into a direct sum of one-dimensional subrepresentations,
but one can verify using \cite{W2} that this does not happen
for either of the 2 two-dimensional irreducible representations
of the quantum group $C(S_3) \rtimes_{\tau} {\Bbb Z}/ 2{\Bbb Z}$.
Similarly, we see that the group of one-dimensional representations
of the algebra $L(G)$ of \cite{Nik} (see Example 2.9.(iii) therein) is
${\Bbb Z} / 4{\Bbb Z}$, while using \cite{W2} we see that
the group of one dimensional representations of the quantum group
$C(S_3) \rtimes_{\tau} {\Bbb Z} / 2{\Bbb Z}$ is
isomorphic to ${\Bbb Z}/ 2{\Bbb Z} \oplus {\Bbb Z}/ 2{\Bbb Z}$.

\section{
Deformation of finite groups of Lie type}

Instead of constructing $q$-deformations of finite groups of Lie type
(see \probref{finitegroupsoflietype}), we now construct a
Rieffel type deformation for such finite groups.

Let $q=p^d$ be the power of a prime $p$,
and let $F_q$ be the field with $q$ elements.
Let $G$ be a finite groups of Lie type over $F_q$, e.g.,
$GL(n, F_q)$, $SL(n, F_q)$, $PSL(n, F_q)$ or $Sp(n, F_q)$.
Let $T$ be a maximal abelian subgroup of $G$ ($T$ can be
viewed as the points over $F_q$ of a maximal torus of
the algebraic group corresponding to $G$).
Apart from trivial cases like $GL(2, F_2) = SL(2, F_2) = PSL(2, F_2) = S_3$,
$T$ has non-trivial skew-symmetric automorphisms $S$.
Letting $J = S \oplus (-S)$, we can use
\thref{theoremquantumgroup} to obtain a deformation
$C(G)_J$. In general, the quantum groups $C(G)_J$ should
be non-trivial in the sense of the previous sections.

Take for example the simplest case $GL(2, F_q)$. Then
$T$ consists matrices of the form
$$  \left(
\begin{array}{cc}
a & 0    \\
0 & b
\end{array}
\right)
 $$
with $a , b \in F^*_q$. Hence $T$ is isomorphic
to the group $({\Bbb Z}/(q-1){\Bbb Z})^2$.
As in the previous section, we have on $T$ the canonical skew-symmetric
automorphism
$$ S = \left(
\begin{array}{cccc}
0 & 1    \\
-1 & 0
\end{array}
\right).
$$
Using \thref{theoremquantumgroup}, we can form the deformation
$C(GL(2, F_q))_J$. If $q \neq 2, 3$, then the remark
after \propref{specsubspace} shows that
$C(GL(2, F_q))_J$ is noncommutative and noncocommutative.

Apart from trivial cases such as those mentioned in the beginning of this
section, it is not hard to see that $C(G)_J$ is noncommutative and
noncocommutative for an arbitrary finite group of Lie type $G$.
However, it is not so clear how to show that $C(G)_J$ cannot
be expressed as crossed product of the form $C(G_1) \rtimes_\tau G_2$
for some finite groups $G_1$ and $G_2$.
We plan to study this further in the future.

Since there are many skew-symmetric automorphisms on
the maximal abelian subgroup of a finite group of Lie type $G$,
$C(G)_J$ is not the analog of $q$-deformation of
Drinfeld and Jimbo of infinite Lie groups.
\vv
{\bf Acknowlegdement.}
The author would like to thank Professor Marc Rieffel for
continual encouragement and support, who urged the author
to clarify the relationship of the construction in this paper
and the construction of Enock and Vainerman \cite{EnockVain1}.
He also thanks C.-S. Chu, A. Weinstein
and S.L. Woronowicz for their helpful comments.
The main results of this paper were obtained while the author was
visiting IHES during the year July, 1995-Aug, 1996. He
thanks the IHES for its financial support and hospitality during this period.
The author also wishes to thank the Department of
Mathematics at UC-Berkeley for its support and hospitality
while he held an NSF Postdoctoral Fellowship there
during the final stage of this paper. He thanks the referee
for pointing out an error at the end of the proof of Theorem 3.2
in the original manuscript.


\vspace{0.5cm}
\hfill 
revised, July 22, 1998


\begin{thebibliography}{99}


\bibitem{BS2} Baaj, S. and Skandalis, G.:
{\rm Unitaires multiplicatifs et dualite pour les produits
croises de $C^*$-algebres,}
{\em Ann. Sci. Ec. Norm. Sup.} {\bf 26} (1993), 425-488.

\bibitem{Cn} Connes, A.:
{\em Noncommutative Geometry,}
Academic Press, 1994.

\bibitem{Cn3} Connes, A.:
{\rm Noncommutative geometry and reality,}
{\em J. Math. Phys.} {\bf 36}:11 (1995), 6194-6231

%

\bibitem{Dr6} Drinfeld, V. G. :
{\rm Quasi-Hopf algebras,}
{\em Leningrad Math. J.} {\bf 1}(6) (1990), 1419-1457

\bibitem{EnockVain1} Enock, M. and Vainerman, L.:
{\rm Deformation of a Kac algebra by an abelian subgroup,}
{\em Commun. Math. Phys.} {\bf 178} (1996), 571-596.

\bibitem{HR2} Hewitt, E. and Ross, K.:
{\em Abstract Harmonic Analysis}   II,
Springer--Verlag,      1970.

\bibitem{Kac3} Kac, G.:
{\rm Certain arithmetic properties of ring groups,}
{\em Funct. Anal. Appl.} {\bf 6} (1972), 158-160.

\bibitem{KP2} Kac, G. and Palyutkin, V.:
{\rm Finite ring groups,}
{\em Trans. Moscow Math. Soc.} {\bf 15} (1966), 251-294.

\bibitem{Landstad1} Landstad, M.B.:
{\rm Quantizations arising from abelian subgroups,}
{\em Internat. J. Math.} {\bf 5} (1994), 897-936.

\bibitem{LandstadRaeburn1} Lanstad, M. B. and Raeburn, I.:
{\rm Twisted dual-group algebras: Equivariant deformations of $C_0(G)$},
{\em J. Funct. Anal.} {\bf 132}, 43-85 (1995).

\bibitem{Lev1} Levendorskii, S.:
{\rm Twisted algebra of functions on compact quantum
group and their representations,}
{\em St. Petersburg Math. J.} {\bf 3}:2 (1992), 405-423.

\bibitem{LS1} Levendorskii, S. and Soibelman, Y.:
{\rm Algebra of functions on compact quantum groups,
Schubert cells, and quantum tori,}
{\em Commun. Math. Phys.} {\bf 139}, (1991), 141-170.

\bibitem{Neumann1} von Neumann, J.:
{\rm Die Eindeutigkeit der Schr\"odingerschen Operatoren,}
{\em Math. Ann.} {\bf 104} (1931), 570-578.

\bibitem{Nik} Nikshych, Dmitri:
{\rm $K_0$ rings and twisting of finite dimensional
semisimple Hopf algebras,}
{\em Preprint}, National Technical University of Ukraine
``Kiev Polytechnic Institute'', 1997.

\bibitem{Rf6} Rieffel, M.:
Deformation quantization for actions of ${\mbox{\bf R}}^d$,
{\em Memoirs A.M.S.} no. {\bf 506}, 1993.

\bibitem{Rf8} Rieffel, M.:
{\rm Compact quantum groups associated with toral subgroups,}
{\em Contemp. Math.} {\bf 145} (1993), 465-491.

\bibitem{Rf9} Rieffel, M.:
{\rm Non-compact quantum groups associated with abelian subgroups,}
{\em Commun. Math. Phys.} {\bf 171} (1995), 181-201.

\bibitem{Vainerman2} Vainerman, L. I.:
{\rm $2$-cocycles and twisting of Kac algebras},
{\em Commun. Math. Phys.} {\bf 191}:3 (1998), 697-721.

\bibitem{VS} Vaksman, L. and Soibelman, Y.:
{\rm The algebra of functions on quantum $SU(2)$,}
{\em Funct. Anal. ego Pril.} {\bf  223} (1988), 1-14.

\bibitem{W5} Van Daele, A. and Wang, S. Z.:
{\rm Universal quantum groups,}
{\em International J. Math} {\bf 7}:2 (1996), 255-264.

\bibitem{W1} Wang, S. Z.:
{\rm Free products of compact quantum groups,}
{\em Commun. Math. Phys.} {\bf 167}:3 (1995), 671-692.

\bibitem{W2} Wang, S. Z.:
{\rm Tensor products and crossed products of compact quantum groups,}
{\em Proc. London Math. Soc.} {\bf 71}:3 (1995), 695-720.

\bibitem{W3} Wang, S. Z.:
{\rm Krein duality for compact quantum groups,}
{\em J. Math. Phys.} {\bf 38}:1 (1997), 524-534

\bibitem{W4} Wang, S. Z.:
{\rm Deformations of compact quantum groups via Rieffel's quantization,}
{\em Commun. Math. Phys.} {\bf 178}:3 (1996), 747-764.

%

\bibitem{W7} Wang, S. Z.:
{\rm Problems in the theory of quantum groups,}
in {\em Quantum Groups and Quantum Spaces},
Banach Center Publication 40 (1997),
Inst. of Math., Polish Acad. Sci.,
Editors: R. Budzynski, W. Pusz, and S. Zakrzewski.
pp67-78

%

\bibitem{Wor4} Woronowicz, S. L.:
{\rm Twisted $SU(2)$ group. An example of noncommutative
differential calculus,}
{\em Publ. RIMS, Kyoto Univ.} {\bf 23} (1987), 117-181.

\bibitem{Wor5} Woronowicz, S. L.:
{\rm Compact matrix pseudogroups,}
{\em Commun. Math. Phys.} {\bf 111} (1987), 613-665.

\bibitem{Wor6} Woronowicz, S. L.:
{\rm Tannaka-Krein duality for compact matrix pseudogroups.
Twisted $SU(N)$ groups,}
{\em Invent. Math.} {\bf 93} (1988), 35-76.


\end{thebibliography}
\end{document}